\documentclass[a4paper, 11pt]{article}
\usepackage[margin=0.9in]{geometry}
\usepackage{times}
\usepackage{amsmath,amssymb,amsthm}

\usepackage[utf8x]{inputenc}
\usepackage{epsfig}
\usepackage{graphicx,subfigure}
\usepackage{hyperref,url}
\hypersetup{colorlinks=true,citecolor=black,filecolor=black,linkcolor=black,urlcolor=black} 
\usepackage[ruled,vlined]{algorithm2e}
\usepackage[usenames,dvipsnames]{color}

\renewcommand{\CommentSty}{\hypersetup{linkcolor=Gray} \textcolor{Gray}
  \# \textcolor{Gray}} 

\newcommand{\m}{\mathcal}

\newcommand{\DD}{\displaystyle D \hspace{-0.5mm} D}

\newcommand{\Dl}{\displaystyle D \hspace{-0.2mm} l}
\newcommand{\bs}{\displaystyle \boldsymbol}

\newcommand{\cp}{\textrm{cp}}

\begin{document}
\title{Adaptive Wavelet Collocation Method\\ for Simulation of
  Time Dependent Maxwell's Equations}
\author{Haojun~Li$^*$, Kirankumar~R.~Hiremath$^{*,+}$,
  Andreas~Rieder$^{*}$, and Wolfgang~Freude$^\star$\\
\vspace{0.05cm}\\
$*$Department of Mathematics, Karlsruhe Institute of
  Technology, Germany\\
$+$ Computational Nanooptics Group, Department of Numerical Analysis and Modelling\\ Konrad-Zuse-Zentrum f\"ur
  Informationstechnik Berlin, Takustrasse 7, 14195 Berlin, Germany\\
$\star$Institute of High-Frequency and Quantum Electronics,\\
Karlsruhe Institute of Technology, Germany \\
\vspace{0.05cm}\\
Corresponding author: andreas.rieder@.kit.edu 
}
\date{5 April, 2012}
\maketitle

\begin{abstract}
  This paper investigates an adaptive wavelet collocation time domain
  method for the numerical solution of Maxwell's equations. In this method a
  computational grid is dynamically adapted at each time step by using the
  wavelet decomposition of the field at that time instant. In the regions
  where the fields are highly localized, the method assigns more grid
  points; and in the regions where the fields are sparse, there will be
  less grid points. On the adapted grid, update schemes with high spatial order and explicit
  time stepping are formulated. The
  method has high compression rate, which substantially reduces the
  computational cost allowing efficient use of computational
  resources. This adaptive wavelet collocation method is especially
  suitable for simulation of guided-wave optical devices.
\end{abstract}

\noindent \textbf{keyword}: Maxwell's equations, time domain methods, wavelets,
wavelet collocation method, adaptivity

\section{Introduction}
The numerical solution of Maxwell's equations is an active area of
computational research. Typically, Maxwell's equations are solved
either in the frequency domain or in the time domain, where each of these
approaches has its own relative merits. We are specifically interested
in efficient algorithms for light propagation problems in guided wave
photonic applications~\cite{okamoto_00}, and work in the time 
domain. The most popular class of methods in this area is the finite
difference time domain (FDTD) method~\cite{Taflove05}. Due to the structured
grid requirement of these methods, they become cumbersome while dealing with
optical devices having curved interfaces and different length scales. To
overcome these difficulties, a discontinuous Galerkin time domain (DGTD)
method has been investigated~\cite{hesthaven_08}. For a time dependent
wave propagation problem, all these methods use a fixed grid/mesh for
discretization. In general, such a grid can under-sample the temporal
dynamics, or over-sample the field
propagation causing high  computational costs. If the spatial grid
adapts itself according to the temporal evolution of the field, then the
computational resources will be used much more efficiently. 

We propose an adaptive-grid method which represents propagating fields at
each time step by a compressed wavelet decomposition, and which
automatically adapts the computational mesh to the changing shape of the
signal. In the initial studies of the wavelet formulation,
the interpolating scaling functions were used for frequency domain waveguide
analysis~\cite{Fujii01}. To the best of our knowledge, the suitability
of the wavelet decompositions for time dependent Maxwell problems has
not  yet been investigated. Vasilyev and his co-authors developed the
adaptive wavelet collocation time domain (AWC-TD) method as a general
scheme to solve evolution equations, and  they successfully verified  the
scheme's effectiveness in the area of computational fluid
dynamics~\cite{Vasilyev00, Vasilyev03}. Based on these studies,  we
present in this work a proof-of-concept for
an AWC-TD for the time dependent Maxwell's equations.

The paper is organized as follows. In Sec.\,\ref{sec:TDME}, we provide a
brief account on Maxwell's equations and some of the related concepts
for their numerical solutions. We start Sec.\,\ref{sec:awcm} with an introduction to
(interpolating) wavelets, and how they can be used to discretize partial
differential equations. Also in this section we explain the structure of
AWC-TD method in the context of Maxwell's
equations. Sec.\,\ref{sec:algo} gives algorithmic details of the
method. Numerical results of the AWC-TD method are given in
Sec.\,\ref{sec:ns} which  contains our  
numerical experiments of propagating a 2D Gaussian peak in homogeneous
environment. Finally we close the paper with concluding remarks in
Sec.\,\ref{sec:con}.

\section{Time domain Maxwell's equations}
\label{sec:TDME}
Propagation of optical waves in a linear, non-magnetic dielectric medium
with no charges and currents  is governed by the following time
dependent Maxwell's equations 
\begin{equation}
\label{eq:tMax}
- \dfrac{\partial }{\partial t}\vec{\m{B}}(\vec{r}, t)  =  \nabla \times
\vec{\m{E}}(\vec{r}, t),\quad\dfrac{\partial}{\partial t}
\vec{\m{D}}(\vec{r}, t)  =  \nabla \times \vec{\m{H}}(\vec{r},
t),\quad\nabla \cdot \vec{\m{D}}(\vec{r}, t)  = 0,\ 
\text{ and }\ \nabla \cdot \vec{\m{B}}(\vec{r}, t) = 0,  
\end{equation}
where the electric field $\vec{\m{E}}$ and the electric flux density
$\vec{\m{D}}$,
as well as the magnetic field $\vec{\m{H}}$ and the magnetic flux density
$\vec{\m{B}}$, are related by the constitutive relations
\begin{equation*}
\vec{\m{D}}(\vec{r}, t)= \varepsilon_0 \varepsilon_r(\vec{r}) \vec{\m{E}}(\vec{r}, t)\ \text{ and }\ \vec{\m{B}}(\vec{r}, t) = \mu_0 \vec{\m{H}}(\vec{r}, t).
\end{equation*}
Here $\varepsilon_0$ is the free space permittivity, $\varepsilon_r$ is
the relative permittivity and $\mu_0$ is free space permeability.

For illustration purpose, we restrict ourselves to a 2D setting where
the fields  and the material properties are assumed to be invariant in
the $y$-direction, i.\,e. $\vec{r} = (x, z)$ and the partial derivatives of
all fields with respect to $y$ vanish identically. We suppress the
explicit function dependence on $\vec{r}$ and $t$. Then Maxwell's
equations~\eqref{eq:tMax} decouple into a pair of independent sets of
equations, 
\begin{equation}
  \label{eq:TEy}
\frac{\partial \m{E}_{x}}{\partial t}  = -\frac{1}{\varepsilon_0 \varepsilon_r}
\frac{\partial \m{H}_{y}}{\partial z},~\hfill 
\frac{\partial \m{E}_{z}}{\partial t} = \frac{1}{\varepsilon_0 \varepsilon_r} \frac{\partial \m{H}_{y}}{\partial x},~\hfill
\frac{\partial \m{H}_{y}}{\partial t} = \frac{1}{\mu_{0}}\left(\frac{\partial \m{E}_{x}}{\partial z}-\frac{\partial \m{E}_{z}}{\partial x} \right),  
\end{equation}
identified as \textit{transverse electric} (TE)$_y$ setting, and
\begin{equation}
  \label{eq:TMy}
\frac{\partial \m{H}_x}{\partial t}  = \frac{1}{\mu_{0}} \frac{\partial \m{E}_y} {\partial z}, ~\hfill 
\quad \frac{\partial \m{H}_z}{\partial t}  = -\frac{1}{\mu_{0}} \frac{\partial \m{E}_y} {\partial x},~\hfill  
\frac{\partial \m{E}_y}{\partial t}  = \frac{1}{\varepsilon_0 \varepsilon_r} \left(\frac{\partial  \m{H}_x} {\partial z}-\frac{\partial \m{H}_z}{\partial x} \right),  
\end{equation}
identified as \textit{transverse magnetic} (TM)$_y$ setting. Here $
\m{E}_{x}$, $\m{E}_{z}$, $\cdots$ etc.\ denote the respective field
components. 

\begin{figure}[ht]
  \centering
  \includegraphics[width=10cm]{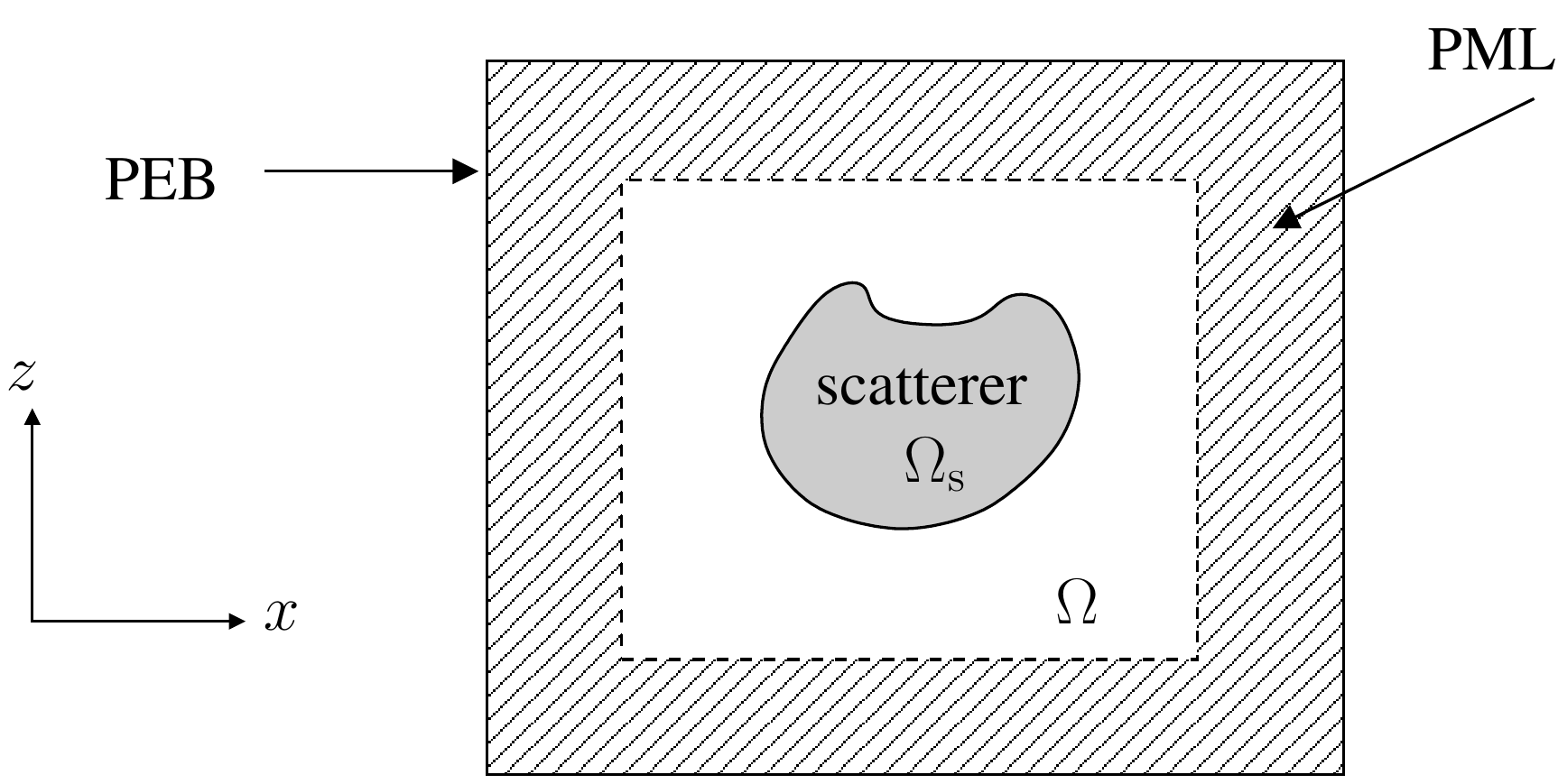}
  \caption{\label{fig:comp_domain} Typical simulation setting with a
    computational domain $\Omega$ surrounded by the perfectly matched
  layer. Here just for the sake of illustration, we 
  show the scatterer $\Omega_{\rm{s}}$ completely enclosed inside $\Omega$. Other
  configurations like incoming-outgoing waveguides  are also
  possible~\cite{HaojunLi_thesis}.}
\end{figure}

Originally,  Maxwell's equations are formulated for a whole space. For
numerical computations we need to restrict them to a bounded
computational domain $\Omega$ as shown in
Fig.~\ref{fig:comp_domain}. This is done with a transparent boundary
condition, which is realized in our case  with {\it perfectly matched
  layer} (PML)~\cite{Berenger94, Gedney96}. The principle of PML is that 
(outgoing) waves scattered from the scatterer $\Omega_{\mathrm{s}}$ pass
through the interface between $\Omega$ and PML without reflections, 
and attenuate significantly inside the PML. The waves virtually vanish before
reaching the outermost boundary of the PML, where the perfectly electric
boundary  (PEB) condition is employed. Implementation details about the PML 
technique specific for the method discussed in this paper can be found
in Ref.\,\cite{HaojunLi_thesis}. For the sake of clarity, we work with
the general formulation given by Eq.\,\eqref{eq:TEy}-\eqref{eq:TMy}.

As in the case of the standard FDTD method~\cite{Taflove05}, in our
approach we use the central difference scheme for the time derivatives in 
Eq.\,\eqref{eq:TEy}-\eqref{eq:TMy}, but we will construct a different
discretization scheme of the spatial derivatives. This is done with
interpolating scaling functions and lifted interpolating wavelets
(explained in Sec.\,\ref{sec:awcm}). The induced multiresolution
approximation~\cite{Mallat89,Mallat98} enables us to 
decompose fields into various resolution levels, and thus allows to discard
unimportant features. As a result, we will obtain a variant of the FDTD
method, which is constructed with respect to a locally refined
grid. In the next section we describe this numerical scheme in detail.

\section{Adaptive wavelet collocation method}
\label{sec:awcm}
The adaptive wavelet collocation (AWC) method was proposed by Vasilyev
and co-authors in a series of papers~\cite{Vasilyev96, Vasilyev97,
  Vasilyev00, Vasilyev03} as a general scheme to solve evolution
equations. In the present section, we tailor the AWC method to tackle
Eq.\eqref{eq:TEy}-\eqref{eq:TMy}. In contrast to the 
originally formulated AWC method, we do not need to utilize  second
generation wavelets, which have been mainly invented to implement 
boundary constraints, and to find wavelet decompositions on irregular 
domains. Since we use the PML method, we can identify  field values
outside the PML region with zero, and therefore we are not forced to adapt our
wavelets to the boundary restrictions. Hence, we consider only the first
generation wavelets, which are generated by the shifts and the dilations of a
single function. Now we outline the essential steps for computing
spatial derivatives of functions in  wavelet representations.

\subsection{Preliminaries}
A starting point of the AWC method is a wavelet decomposition of a
function $f\in L^2(\mathbb{R})$: 
\begin{equation}
\label{eq:wavelet_decomp}
f = \sum_{k \in \mathbb{Z}} \alpha_{j_0,k} \phi_{j_0,k} + \sum_{j=j_0}^{+\infty} 
\sum_{m \in \mathbb{Z}} \beta_{j,m} \psi_{j,m}
\end{equation}
where $j_{0} \in \mathbb{Z}$, $\phi$ is the {\it scaling function} and
$\psi$ is the {\it wavelet function}~\cite{DaubTen92, Rieder98}. For all
$j,n \in \mathbb{Z}$, by $\phi_{j,n}$ and
$\psi_{j,n}$  we abbreviate the dilated and translated versions of $\phi$
and $\psi$, i.e. $\phi_{j,n} (\cdot) = 2^{j/2} \phi(2^j\cdot-n)$,
$\psi_{j,n} (\cdot) = 2^{j/2} \psi(2^j\cdot-n)$.

The first (single) sum in (\ref{eq:wavelet_decomp}) represents \textit{rough} or
  \textit{low frequency information} of  $f$, while the second (double) sum
contains  the \textit{detail information} at various resolution levels
starting from the level $j_0$ to $+\infty$. The absolute magnitude of
the coefficients $\alpha_{j_0,k}$ and $\beta_{j,m}$ measure the
contributions of $\phi_{j_0,k}$ and $\psi_{j,m}$ to $f$. By discarding
terms in the double sum for which the wavelet coefficients $\beta_{j,m}$
are absolutely less than a given threshold, one can efficiently compress
the representation of $f$. This wavelet decomposition compression
principle is exploited in the AWC method to enhance the computational
efficiency.

There are various families of the scaling functions $\phi$ and  wavelet
functions $\psi$  allowing representations like~\eqref{eq:wavelet_decomp}. As
in~\cite{Vasilyev00, Vasilyev03}, we work with the interpolating scaling
functions~\cite{DD89} and the corresponding lifted interpolating
wavelets~\cite{Sweldens96, Sweldens98}. Due to their interpolation
property, we have
\begin{equation*}
\phi(k)=\delta_{0,k}=\begin{cases}
1\ :&\! k=0,\\[1mm]
0\ : &\! k\in\mathbb{Z}\setminus\{0\},
\end{cases}
\end{equation*}
and as a result, there exits a unique grid associated with the family
$\{\phi_{j,k}\}$. The resulting numerical scheme can be seen as a
variant of the well known finite difference method. We exploit this
interpolating property in Sec.\,\ref{sec:wc} and Sec.\,\ref{sec:sp}.

In particular, we use the interpolating scaling function (ISF) family
developed by Deslauriers and Dubuc~\cite{Dubuc86, DD89}. They
constructed the interpolating functions by the {\it iterative
  interpolation method}, which does not require the concept of
wavelets. Later Sweldens~\cite{Sweldens96, Sweldens98} constructed the
corresponding wavelet by lifting the Donoho 
wavelet~\cite{Donoho92}. We use $\DD_N$ to denote ISF of order $N$, and
$\Dl_{\widetilde N}$ to denote the lifted interpolating wavelet of
order~$\widetilde N$. Here the order $N$ means that any polynomial $p$ of degree
$k \leq 2N-1$ can be expressed as 
\begin{equation*} 
  p(\cdot) = \sum_m c_m \DD_N(\cdot - m)
\end{equation*}
with suitable coefficients $\{c_m\}$.
The order $\widetilde{N}$ is half the number of the vanishing moments of the lifted interpolating wavelet, i.e.,
\begin{equation*}
  \int x^k \Dl_{\widetilde{N}} (x) \mathrm{d}x =0, \quad k=0, 1, \dots, 2\widetilde{N}-1.
\end{equation*}
Further details can be found in \cite{Sweldens96, Sweldens98,
  HaojunLi_thesis}. We 
normally choose same orders for the ISF and the lifted interpolating
wavelet, i.e., $N={\widetilde N}$. It is easy to see that $\DD_N$ and
$\Dl_{\widetilde{N}}$ have compact supports, which increase with
the order $N$.

For the TM$_y$ setting in Eq.~\eqref{eq:TMy}, the electric  and magnetic
fields depend on the spatial variables $(x,z)$. As usual, see, e.g., \cite{Mallat98,Rieder98}, we represent 2D fields by expansions of 2D scaling functions and wavelets which are defined by
\begin{align*}
  {\bs \phi}_N (x,z) & := \DD_N (x) \DD_N (z), \\[2mm]
  {\bs \psi}^{\nu}_N (x,z) & := \left \{\! \begin{array}{ll}
   \Dl_N(x) \DD_N(z) &: \quad \nu = 1,\\[1mm]
   \DD_N(x) \Dl_N(z) &: \quad \nu = 2,\\[1mm]
   \Dl_N(x) \Dl_N(z) &: \quad \nu = 3, \end{array} \right .
\end{align*}
and use the following abbreviations
\begin{align*}
  ({\bs \phi}_N)_{j,m,n} (x,z) & := {(\DD_N)}_{j,m} (x) {(\DD_N)}_{j,n} (z), \\[2mm]
  ({\bs \psi}_N^{\nu})_{j,m,n} (x,z) & := \left\{\!  \begin{array}{ll}
   {(\Dl_N)}_{j,m} (x) {(\DD_N)}_{j+1,2n} (z) &: \quad \nu = 1, \\[1mm]
  {(\DD_N)}_{j+1,2m} (x) {(\Dl_N)}_{j,n} (z) &: \quad \nu = 2, \\[1mm]
  {(\Dl_N)}_{j,m} (x) {(\Dl_N)}_{j,n} (z) &: \quad \nu = 3. \end{array} \right.
\end{align*}
Let $j_{\min}$ and $j_{\max}$ (with $j_{\min} < j_{\max}$) be the
coarsest and the finest spatial resolution levels. 
Let us consider $f\in L^2(\mathbb{R}^2)$ with exact resolution level $j_{\max}$, that is,
\begin{equation}
\label{eq:2d_scaling_rep}
  f = \sum_{m,n} \alpha_{j_{\max},m,n} ({\bs \phi}_N)_{j_{\max},m,n}.
\end{equation}
Then the wavelet representation of $f$ with coarsest resolution level $j_{\min}$ is given by
\begin{equation}
\label{eq:2d_decomposition}
  f = \sum_{m,n} \alpha_{j_{\min},m,n} ({\bs \phi}_N)_{j_{\min},m,n} +
  \sum_{\nu = 1}^{3} \sum_{j = j_{\min}}^{j_{\max} - 1} \sum_{m,n}
  \beta^{\nu}_{j,m,n} ({\bs \psi}_N^{\nu})_{j,m,n}
\end{equation}
where the scaling coefficients $\{\alpha_{j_{\min},m,n}\}$ and the wavelet
coefficients $\{\beta^{\nu}_{j,m,n}\}$ can be calculated from the level
$j_{\max}$ scaling coefficients $\{\alpha_{j_{\max},m,n}\}$ by the
\textit{normalized} 2D \textit{forward wavelet transform} (FWT):
 \begin{subequations}
\label{eq:2dfwt_d}
   \begin{align}
    d^1_{j,m,n} &=  \frac{1}{2} \Big( c_{j+1,2m+1,2n} - \sum_l 2 \tilde{s}_{-l} c_{j+1,2m+2l,2n}
 \Big),\label{eq:2dfwt_d1}\\
    d^2_{j,m,n} &=  \frac{1}{2} \Big( c_{j+1,2m,2n+1} - \sum_l 2 \tilde{s}_{-l} c_{j+1,2m,2n+2l}
 \Big),\label{eq:2dfwt_d2}\\
    d^3_{j,m,n} &=  \frac{1}{4} \Big( c_{j+1,2m+1,2n+1} - \sum_l 2 \tilde{s}_{-l} c_{j+1,2m+2l,2n+1} 
- \sum_{l'} 2 \tilde{s}_{-l'} c_{j+1,2m+1,2n+2l'} \nonumber \\
& + \sum_l \sum_{l'} (2 \tilde{s}_{-l})(2 \tilde{s}_{-l'}) c_{j+1,2m+2l,2n+2l'} \Big),\label{eq:2dfwt_d3}\\
    c_{j,m,n} &=  c_{j+1,2m,2n} + \sum_l s_{-l} d^1_{j,m+l,n} + \sum_{l'} s_{-l'} d^2_{j,m,n+l'} + \sum_l
 \sum_{l'} s_{-l} s_{-l'} d^3_{j,m+l,n+l'},\label{eq:2dfwt_c}
  \end{align}\label{eq:2dfwt}
\end{subequations}
with the following normalization conventions
$$
c_{j,m,n} = 2^j \alpha_{j,m,n}, \quad d^1_{j,m,n} = 2^{j+1/2} \beta^1_{j,m,n}, \quad d^2_{j,m,n} =  2^{j+1/2} \beta^2_{j,m,n} \text{ and } d^3_{j,m,n} = 2^j \beta^3_{j,m,n}. 
$$
The coefficients $2\tilde{s}_{l}$ and $s_{l}$ are Lagrangian interpolation weights. For example, when $N=2$,
these weights are 
$$
s_{-2}=-1/16,\quad s_{-1}=9/16,\quad s_0=9/16,\quad s_1=-1/16, \quad2\tilde{s}_{-1}=-1/16,
$$
and
$$
 2\tilde{s}_0=9/16,\quad 2\tilde{s}_1=9/16, \quad 2\tilde{s}_2=-1/16. 
$$
Readers may consult~\cite{DD89,Goedecker} and \cite[Theorem~12]{Sweldens96} 
for an explanation of how and why Lagrangian weights enter the iterative
interpolation process. 

We also can compute back from the wavelet
representation~\eqref{eq:2d_decomposition} to the scaling function
representation~\eqref{eq:2d_scaling_rep} by the {\em inverse wavelet
  transform} (IWT): 
 \begin{subequations}
   \begin{align}
    c_{j+1,2m,2n} &=  c_{j,m,n} - \sum_l s_{-l} d^1_{j,m+l,n} + \sum_{l'} s_{-l'} d^2_{j,m,n+l'} + \sum_l \sum_{l'} s_{-l} s_{-l'} d^3_{j,m+l,n+l'},\label{eq:2diwt_c}\\
    c_{j+1,2m+1,2n}& =  2d^1_{j,m,n} + \sum_l 2 \tilde{s}_{-l} c_{j+1,2m+2l,2n}, \label{eq:2diwt_d1}\\
    c_{j+1,2m,2n+1} &=  2d^2_{j,m,n} + \sum_l 2 \tilde{s}_{-l} c_{j+1,2m,2n+2l}, \label{eq:2diwt_d2}\\
    c_{j+1,2m+1,2n+1}& =  4d^3_{j,m,n} + \sum_l 2 \tilde{s}_{-l}
    c_{j+1,2m+2l,2n+1} + \sum_{l'} 2 \tilde{s}_{-l'} c_{j+1,2m+1,2n+2l'}
    \nonumber \\
& - \sum_l \sum_{l'} (2 \tilde{s}_{-l})(2 \tilde{s}_{-l'}) c_{j+1,2m+2l,2n+2l'}. \label{eq:2diwt_d3}
  \end{align}
  \end{subequations}

\subsection{Adaptive grid refinement  wavelet compression}
\label{sec:wc}
We thin out the triple sum in (\ref{eq:2d_decomposition}) by discarding
small wavelet coefficients, which corresponds to small scale
details. For a given threshold $\zeta >0$, let 
\begin{equation*}
  f_{\zeta} := \sum_{m,n} \alpha_{j_{\min},m,n} ({\bs
    \phi}_N)_{j_{\min},m,n} + \sum_{\nu = 1}^{3} \sum_{j =
    j_{\min}}^{j_{\max} - 1} \sum_{m,n}
  T^{\nu}_{\zeta}(\beta^{\nu}_{j,m,n}) ({\bs \psi}_N^{\nu})_{j,m,n}, 
\end{equation*}
where the threshold function
$T^{\nu}_{\zeta}\colon\mathbb{R}\to\mathbb{R}$ is defined by 
\begin{equation*}
 T^{\nu}_{\zeta}(x) = \left \{\! 
\begin{array}{rl} 
x &: \text{ for } \nu \in\{ 1,2\} \mbox{ and } \ |x| \geq 2^{-j - 1/2}
\zeta,\\
x &:  \text{ for } \nu=3 \mbox{ and } \ |x| \geq 2^{-j} \zeta, \\  
0 &:\ \text{otherwise.} \end{array} \right.
\end{equation*} 
Note that we have defined the uniform  threshold $\zeta$ in terms of the
normalized wavelet coefficients $d^{\nu}_{j,m,n}$ defined in
Eq.\,\eqref{eq:2dfwt_d} i.e., if $|d^{\nu}_{j,m,n}| < \zeta$ then 
$d^{\nu}_{j,m,n}=0$ in $f_\zeta$ in Eq.~\eqref{eq:2d_decomposition}.  
Then the compression error is proportional to $\zeta$ \cite{Vasilyev00}:
\begin{equation*}
  \|f-f_{\zeta}\|_{\infty} \leq C \zeta.
\end{equation*}
Our basis functions in \eqref{eq:2d_decomposition}, which are translates
and dilates of ${\bs \phi}_N$ and ${\bs \psi}_N^\nu$, are interpolating
at the corresponding grid points. Let
\begin{equation*}
  x_{j,m}:=\frac{m}{2^j}\ \text{ and }\ z_{j,n}:=\frac{n}{2^j} \text{
    for } m,n \in \mathbb{Z}, 
\end{equation*}
then we have the following one-to-one correspondence between the basis
functions and the grid points: 
\begin{align*}
  &({\bs \phi}_N)_{j,m,n}\longleftrightarrow (x_{j,m},\,z_{j,n}), &
  \quad & ({\bs \psi}^1_N)_{j,m,n}\longleftrightarrow (x_{j+1,2m+1}, \,
  z_{j+1,2n}),\\ 
  &({\bs \psi}^2_N)_{j,m,n}\longleftrightarrow (x_{j+1,2m}, \,
  z_{j+1,2n+1}), & & ({\bs \psi}^3_N)_{j,m,n} 
\longleftrightarrow (x_{j+1,2m+1},\, z_{j+1,2n+1}).
\end{align*}
Here this correspondence means the validity of the interpolation property. For
instance, we have that 
$$({\bs \phi}_N)_{j,m,n}(x_{j,m'},\,z_{j,n'}) =
\delta_{m,m'} \delta_{n,n'}.
$$ 
With this explanations, we justified the synonymous usage of
\emph{compression of the wavelet representation} and
\emph{compression/adaption of the grid points}. 
\subsection{Adjacent zone}
With the above described wavelet compression, the grid gets suitably
sampled only for the current state of the fields. For a meaningful
(i.e. physical) field evolution in the next time-step, the grid need 
to be  supplemented by additional grid points, on which the fields may become
significant in the next time step. This allows the grid to capture
correctly the propagation of a wave. To this end
Vasilyev~\cite{Vasilyev00,Vasilyev03} has introduced a concept of an
{\it adjacent zone}.

\begin{figure}[!h]
  \centering
   \includegraphics[width=10cm]{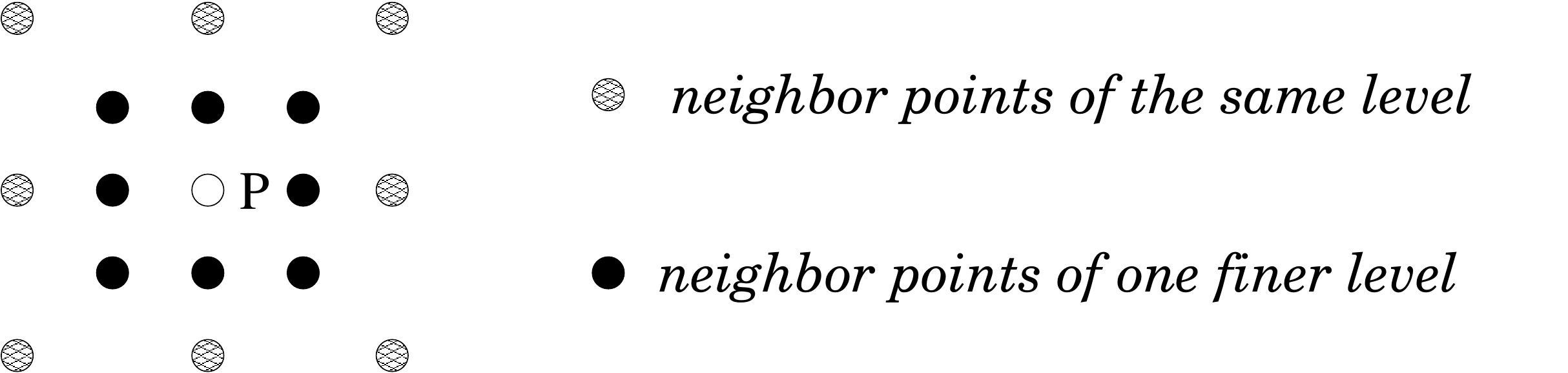}
  \caption{Description of the adjacent zone of a grid point
    $P$.\label{fig:adjacent}}
\end{figure}

To each point $P=(x_{j,m}, z_{j,n})$ in the current grid, 
we attach an adjacent zone which is defined as the set of points
$(x_{j',m'}, z_{j',n'})$ which satisfy 
\begin{equation*}
  |j'-j| \leq L, \quad |2^{j'-j}m-m'| \leq M, \quad |2^{j'-j}n-n'| \leq M,
\end{equation*}
where $L$ is the width of the adjacent levels and $M$ is the width of the
physical space. As in \cite{Vasilyev00}, we verified that $L=M=1$ is 
a computationally sufficient choice. Then the
adjacent zone for a  point $P$ can be depicted as in Fig.~\ref{fig:adjacent}.

Note that the concept of adjacent zone is reasonable only for
continuously propagating waves, as in case of our guided-wave
applications, where in each time step the propagating waves do not
travel far from the current position due to their finite 
propagation speed.

\subsection{Reconstruction check}
\begin{figure}[!ht]
    \centering
    \subfigure[$\times$: The point corresponding to $d^1_{j,m,n}$; \
    $\bullet$: The neighboring points needed to calculate $d^1_{j,m,n}$.]
    {
    \includegraphics[width=4cm]{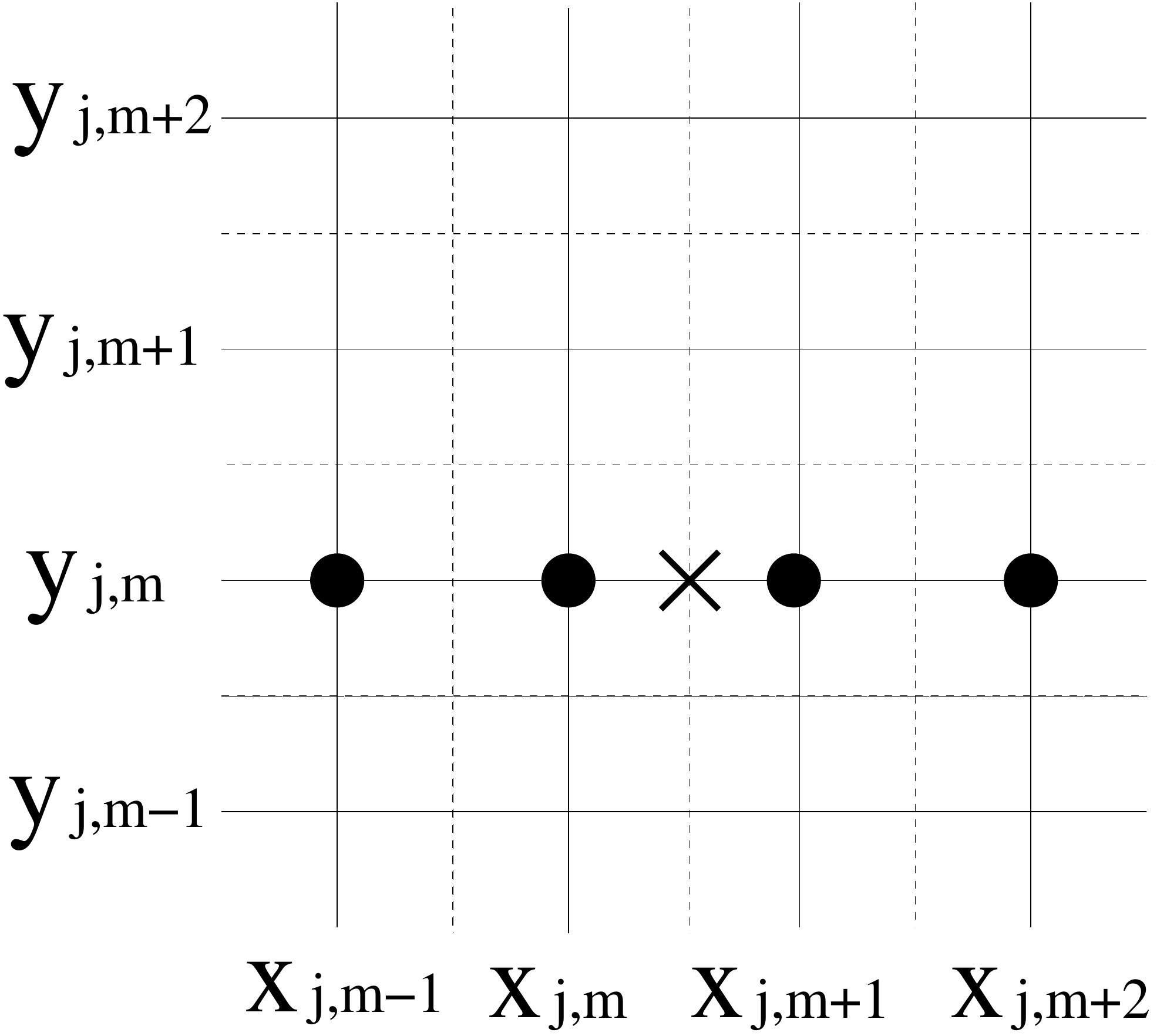}
    \label{fig:grid1}
    }\hspace{0.04\linewidth}
    \subfigure[$\times$: The point corresponding to $d^2_{j,m,n}$; \
    $\bullet$: The neighboring points needed to calculate $d^2_{j,m,n}$.]
    {
    \includegraphics[width=4cm]{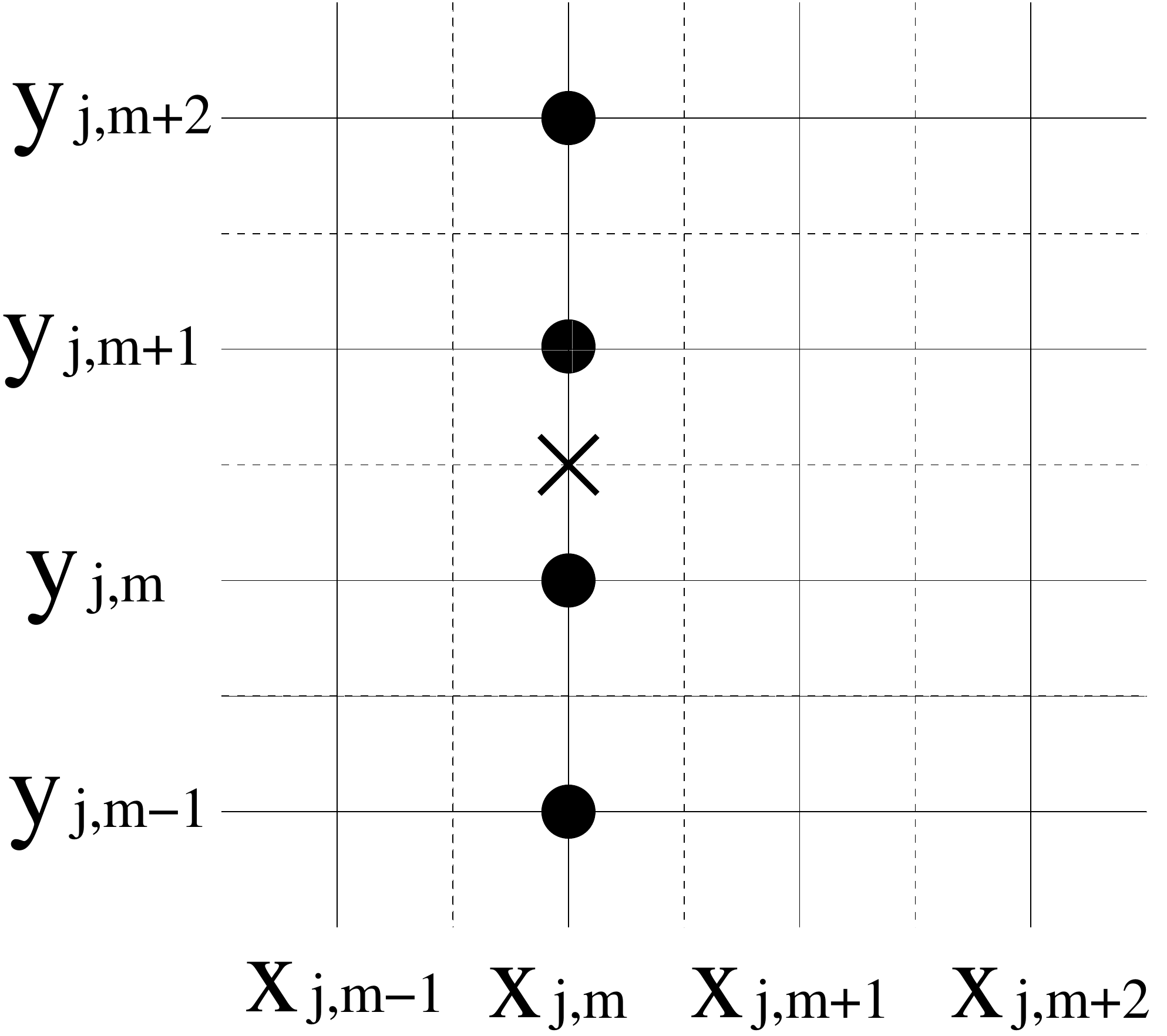}
    \label{fig:grid2}
    }\hspace{0.04\linewidth}
    \subfigure[$\times$: The point corresponding to $d^3_{j,m,n}$; \
    $\bullet$, $\circ$: The neighboring points needed to calculate $d^3_{j,m,n}$.]
    {
    \includegraphics[width=4cm]{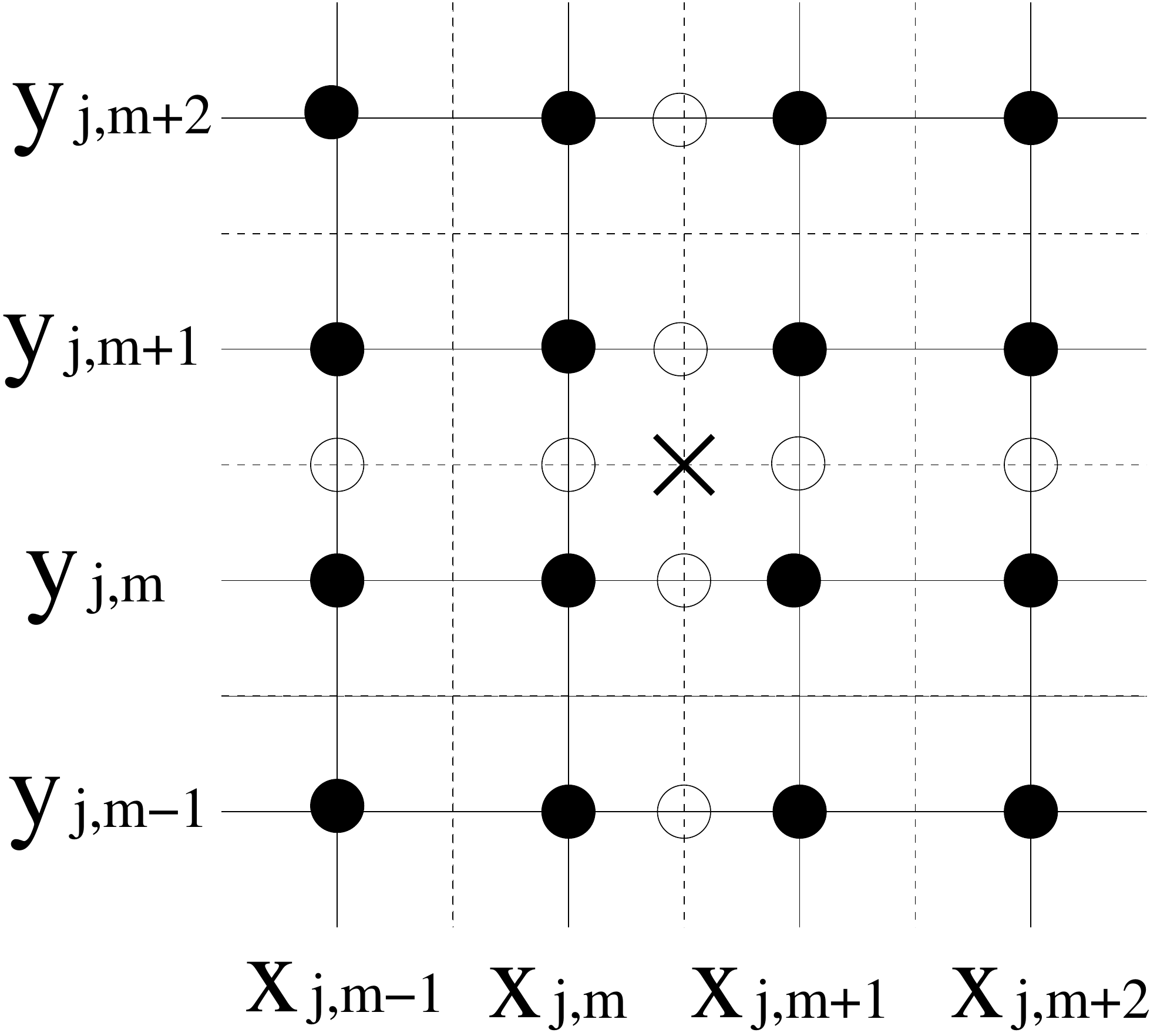}
    \label{fig:grid3}
    }
    \caption{Descriptions of the neighboring points needed to calculate
      the wavelet coefficients $d^{\nu}_{j,m,n}$ with the orders $N=\tilde{N}=2$.}
    \label{fig:wc_grid}
  \end{figure}

In this work we use the wavelet decompositions of the fields only to
determine the adaptive grid. We do not propagate fields in their wavelet
representations (cf.~the statement in the first paragraph of
Sec.~\ref{sec:algo}). Thus at each time step, after adapting the
grid using the FWT, and adding the adjacent zone, we  need to restore the
fields in the physical space by performing the inverse wavelet
transformation (IWT). To this end, we may need to
augment the adaptive grid with additional neighboring points (e.g. see Fig~\ref{fig:wc_grid}). This process of adding neighboring points needed to
calculate the wavelet coefficients in the next time step is called {\it
  reconstruction check}.  Fig~\ref{fig:wc_grid} shows various possible
scenarios, and the corresponding minimal set of the grid points required for
calculation of the wavelet coefficients. The values of the wavelet
coefficients at these {\it newly} added points are set to zero.

The efficiency of the wavelet transform depends on the number
of the finest grid points only at the beginning; however, after the
first compression, it depends solely on the cardinality (= number of
grid points) of the adaptive grid.
\subsection{Calculation of the spatial derivatives on the adaptive grid}
\label{sec:sp}
After the adjacent zone correction and the reconstruction check, we are in a
position to calculate the derivative of $f_{\zeta}$ at a grid point in
the adaptive grid. For this we need to know the {\it density level} of this point, which is defined
as the maximum of the {\it $x$-level} and the {\it $z$-level} of that
point. 

\begin{figure}[!ht]
  \centering
  \includegraphics[width=4cm]{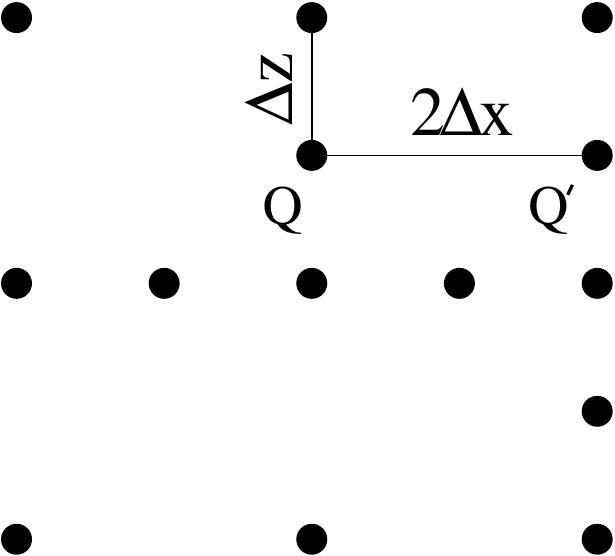}
  \caption{Description of the density level of a point $Q$ in an adaptive
    grid: the $x$-level of $Q$ is $j_{\max}-1$ and the $z$-level of $Q$
    is $j_{\max}$, thus, the density level of $Q$  is
    $j_{\max}$.\label{fig:densitylevel}} 
\end{figure}

We illustrate this concept explicitly only for the $x$-level, the
$z$-level can be determined analogously. For a point $Q = (x_0,z_0)$
in the  adaptive grid $\mathcal{G}$, let $Q' = (x_1, z_0)\in\mathcal{G}$
be the nearest point to $Q$. Then the  $x$-level $Levelx$ of $Q$ relative to $\mathcal{G}$ is  
\begin{equation}\label{eq:xlevel}
  Levelx := j_{\max} - \log_2 (\mathrm{dist}(Q,Q')/\Delta x)
\end{equation}
where $\Delta x$ is the smallest computational mesh size along the $x$
axis, and $\mathrm{dist}(Q,Q') = |x_1 - x_0|$. For $\mathrm{dist}(Q,Q') =
\Delta x$, the level $Levelx$ of $Q$ attains its maximum $j_{\max}$. For
$\mathrm{dist}(Q,Q') = 2 \Delta x$, we have $Levelx=j_{\max}-1$,
etc. See Fig.\,\ref{fig:densitylevel} for an example of describing the 
density level of a grid point. 
  
Now we continue to discuss the derivative calculations. Suppose $j_0$ to
be the density level of $Q$ in $\mathcal{G}$. Then, we can represent
$f_{\zeta}$ by a finite sum ${\bf P}_{j_0} f$ locally in some
neighborhood $\Omega_0$ of $Q$. 
\begin{equation} \label{eq:awcm2d_diff_appr}
  {\bf P}_{j_0} f (x,z) = \sum_{m,n} \alpha_{j_0,m,n} ({\bs
    \phi}_N)_{j_0,m,n} (x,z), \quad (x,z) \in \Omega_0 
\end{equation}
We differentiate ${\bf P}_{j_0} f$ with respect to $x$ to approximate
the $x$-derivative of $f$ at $Q$. If any points in the sum
(\ref{eq:awcm2d_diff_appr}) are not present in $\mathcal{G}$, then we
interpolate the values at these points by the IWT using the values of
the coarser levels. From the interpolation property of $({\bs
  \phi}_N)_{j_0,m,n}$ we know that 
\begin{equation*}
\alpha_{j_0,m,n} = 2^{-j_0} ({\bf P}_{j_0} f) \Big( \frac{m}{2^{j_0}},
\frac{n}{2^{j_0}} \Big), \quad \text{for } m,n \in \mathbb{Z}. 
\end{equation*}
Thus, we have
\begin{equation} \label{eq:awcm_diff}
({\bf P}_{j_0} f) (x,z) = \sum_{m,n} ({\bf P}_{j_0} f) \Big(
\frac{m}{2^{j_0}}, \frac{n}{2^{j_0}} \Big) \DD_N \big( 2^{j_0}x - m
\big) \DD_N \big( 2^{j_0}z - n \big), \quad (x,z) \in \Omega_0. 
  \end{equation}
Differentiate both sides of (\ref{eq:awcm_diff}) with respect to $x$
gives
  \begin{equation} \label{eq:awcm2d_diff1}
\frac{\partial ({\bf P}_{j_0} f)}{\partial x} (x,z) = \sum_{m,n} ({\bf
  P}_{j_0} f) \Big( \frac{m}{2^{j_0}}, \frac{n}{2^{j_0}} \Big)
\frac{\mathrm{d} \DD_N \big( 2^{j_0}x - m \big)}{\mathrm{d}x} \DD_N
\big( 2^{j_0}z - n \big), \quad (x,z) \in \Omega_0. 
\end{equation}
The derivatives of $\DD_N'$ can be calculated exactly at the integers using
the difference filters shown in Table~\ref{table:deri_filter} (see
Ref.~\cite{DD89} for details of the derivation).

  \begin{table}[!h]
    \begin{center}
    \begin{tabular}{|c|r|r|r|} 
    \hline
    $i$  & $N=2$   & $N=3$      & $N=4$           \\
    \hline
    $1$  & $2/3$   &  $272/365$ & $39296/49553$     \\
    $2$  & $-1/12$ &  $-53/365$ & $-76113/396424$  \\
    $3$  &         & $16/1095$  & $1664/49553$       \\
    $4$  &         & $1/2920$   & $-2645/1189272$   \\
    $5$  &         &            & $-128/743295$     \\
    $6$  &         &            & $1/1189272$ \\
    \hline
    \end{tabular}
\caption{Difference filters $\{\DD_N'(i)\}_{i\in\mathbb{Z}}$ with consistency order $2N$. Note that $\DD_N'(i)=-\DD_N'(-i)$. 
\label{table:deri_filter} }
    \end{center} 
  \end{table}

Since the density level of $Q$ is $j_0$, there exist $m',n'
\in \mathbb{Z}$ such that $Q=\big(\frac{m'}{2^{j_0}}, \frac{n'}{2^{j_0}}
\big)$ and it is easy to see that  
\begin{align}
\frac{\partial ({\bf P}_{j_0} f)}{\partial x}
\Big(\frac{m'}{2^{j_0}},\frac{n'}{2^{j_0}} \Big) &= \sum_{m,n} ({\bf
  P}_{j_0} f) \Big( \frac{m}{2^{j_0}}, \frac{n}{2^{j_0}} \Big)
\frac{\mathrm{d} \DD_N \big( m' - m \big)}{\mathrm{d}x} \DD_N \big( n'
- n \big) \nonumber\\[1mm] 
&=  2^{j_0} \sum_{m} ({\bf P}_{j_0} f) \Big( \frac{m}{2^{j_0}},
\frac{n'}{2^{j_0}} \Big) \DD_N' \big( m'- m \big).\label{eq:diffx} 
\end{align}
Similarly, 
  \begin{equation}
    \frac{\partial ({\bf P}_{j_0} f)}{\partial z}
    \Big(\frac{m'}{2^{j_0}},\frac{n'}{2^{j_0}} \Big) = 2^{j_0}
    \sum_{n} ({\bf P}_{j_0} f) \Big( \frac{m'}{2^{j_0}},
    \frac{n}{2^{j_0}} \Big) \DD_N' \big( n' - n \big). 
    \label{eq:diffz}
  \end{equation}
This finishes the general discussion about the adaptive wavelet
collocation method; in the next section, we apply it to Maxwell's
equations.

\section{AWC-TD method for Maxwell's equations}
\label{sec:algo}
In this section we formulate the update scheme for Maxwell's equations,
and then elaborate on algorithmic issues related with the AWC-TD
method. In the present formulation we represent the 
electric and magnetic fields in the physical space, and not in the
wavelet space. To unleash the full power of adaptivity, however, the
field representation and the update in wavelet space are advantageous. 

We illustrate the method for the transverse magnetic (TM)$_y$ setting given
by \eqref{eq:TMy}. Similar procedure can also be formulated for
TE$_y$ setting in \eqref{eq:TEy}. Unlike the standard FDTD
  method, here the electric field and the magnetic field components are
  evaluated on same spatial grid, and their spatial derivatives are
  approximated at the same grid point. But the electric field components
  are sampled at integer time-steps, whereas the magnetic field
  components are sampled at half-integer time-steps.

\subsection{Update scheme for the spatial derivative}
\label{sec:upsc}

For a point $Q$ in the adapted grid $\mathcal{G}$, let
$\mathcal{H}_x|^{k+1/2}_{Q}$, $\mathcal{H}_z|^{k+1/2}_{Q}$ and
$\mathcal{E}_y|^{k}_{Q}$ denote the discretized value of 
$\mathcal{H}_x$, $\mathcal{H}_z$ and $\mathcal{E}_y$ at the point $Q$, and
at a time $(k+1/2) \Delta t$ for the magnetic field components and at a time
$k \Delta t$ for the electric field component where $\Delta t>0$ is the time step size (Note that, the
electric field components are sampled at integer time-steps, whereas
the magnetic field components are sampled at half-integer
time-steps.). Assume $j(Q)$ to be the density level of $Q$ relative to
$\mathcal{G}$. Then we can represent the point $Q$ as
$(x_{j(Q),m'},z_{j(Q),n'})$ for some $m',n' \in \mathbb{Z}$.

Let $L$ be the length of the computational domain $\Omega$. We rescale
the wavelet decomposition \eqref{eq:awcm_diff} with the factor $L$. Then
using the central difference scheme for the time derivatives and using 
\eqref{eq:diffx}-\eqref{eq:diffz} for the spatial derivatives, we get
the following difference equations
\begin{subequations}\label{eq:update_TMy}
 \begin{align}
\mathcal{H}_x|^{k+\frac{1}{2}}_{Q} &= \mathcal{H}_x|^{k-\frac{1}{2}}_{Q} + \frac{\Delta t}{\mu_{0}}
\frac{2^{j(Q)}}{L} \sum_{n} \mathcal{E}_y|^k_{(x_{j(Q),m'}, z_{j(Q),n})}
\DD_N'(n'-n),\label{eq:update_Hx}\\
\mathcal{H}_z|^{k+\frac{1}{2}}_{Q} &= \mathcal{H}_z|^{k-\frac{1}{2}}_{Q} + \frac{\Delta t}{\mu_{0}}
\frac{2^{j(Q)}}{L} \sum_{m} \mathcal{E}_y|^k_{(x_{j(Q),m}, z_{j(Q),n'})}
\DD_N'(m'-m),\label{eq:update_Hz}\\
\mathcal{E}_y|^{k+1}_{Q} &= \mathcal{E}_y|^{k-1}_{Q} + \frac{\Delta
  t}{\varepsilon_{0}} \frac{1}{\varepsilon_r|_{Q}} 
\frac{2^{j(Q)}}{L} \Big( 
\sum_{n} \mathcal{H}_x|^{k+\frac{1}{2}}_{(x_{j(Q),m'}, z_{j(Q),n})}
\DD_N'(n'-n)  \nonumber \\ 
&  - \sum_{m} \mathcal{H}_z|^{k+\frac{1}{2}}_{(x_{j(Q),m},
    z_{j(Q),n'})} \DD_N'(m'-m) \Big),\label{eq:update_Ey}
\end{align}
\end{subequations}
The first time step ($k=0$) is an explicit Euler step with step size
$\Delta t/2$ using initial conditions for the fields at the time $t=0$.
If not explicitly mentioned, otherwise the fields are set zero at the
beginning for all our numerical  experiments in Sec.\,~\ref{sec:ns}. The
update equations for the PML assisted Maxwell's equations can be found 
in Ref.~\cite{HaojunLi_thesis}. 

From the form of these update equations, it is clear that the AWC-TD method
can be thought as an variant of high order FDTD method. The AWC-TD
method is defined with respect to a locally adapted mesh, and unlike the
FDTD method, it does not require a static (fixed), structured mesh. This
will lead to efficient use of the computational resources. In the
next section, we elaborate on algorithmic aspects of the method.

\subsection{Update scheme for the time derivative}
\label{sec:tupsc}
Several choices are available for time stepping. As
in case of the standard FDTD method, we use in \eqref{eq:update_TMy} the central difference
scheme for the discretization of the time derivatives.  For this
explicit scheme, the smallest spatial step-size restricts the maximal
time-step according to the Courant–Friedrichs–Lewy
(CFL) stability condition.  Using a uniform spatial mesh in the 
update equations \eqref{eq:update_TMy} with a mesh size $\Delta$ 
in both coordinate directions the CFL condition reads
\begin{equation}
\label{eq:CFL}
    \Delta t \leq \frac{\Delta}{\sqrt{2} \mathrm{c} \sum_{l=0}^{l_0-1}
      |\DD_N'(l)|},
  \end{equation} 
see \cite[Sec. 3.5]{HaojunLi_thesis} and \cite{Tentzeris99},
where  $\mathrm{c}$ is the speed of light in vacuum and $\{\DD_N'(l)\}$ is the known derivative
filter of the ISF as in
Table~\ref{table:deri_filter}. Due to the local 
adaptive grid strategy of the AWC-TD method, we cannot define a global
stability criteria as above. But choosing
$\Delta$  to be the smallest step size in the adaptive grid, we get a 
conservative bound for $\Delta t$ via \eqref{eq:CFL} for the AWC-TD method. In the
simulation tests (in this paper, and in \cite{HaojunLi_thesis})
we did not experience any stability related issues with this modus operandi.

\subsection{Implementation aspects}
\subsubsection{Grid management}
\label{sec:gm}
In AWC method the computational grid is changed with the state (spatial
localization) of the propagating field. Thus the grid management is one
of the important steps in the implementation of this method. This is
done as following: We store the information of the adaptive grid into a 2D Boolean array called a {\it grid mask} or simply a {\it
  mask},  whose size is square the number of the finest grid points along one
direction. We use 2D arrays of real numbers with the size of the grid
mask to store the fields such as $\m{E}_y$, $\m{H}_x$, and $\m{H}_z$
etc. Note that the computational effort for updating the fields at each
time step is proportional to the cardinality (i.e. the number of entries
in the mask with value $1$) of the adaptive grid.

If the value of an entry of a mask is $true$ or $1$, then the corresponding
grid point is included in the adaptive grid; otherwise, it is not
included in the grid. Thus by forcing the value of an entry of a
mask to $1$, we  can {\it include} the corresponding point to the grid, or
by forcing the entry  to $0$, we can {\it exclude} the corresponding point
from the grid.

\subsubsection{Algorithmic procedures}
Algorithm~\ref{algo:awcm_TMy} outlines the main function {\it
  awcm\_main()} of AWC-TD method for TM$_y$ setting. It mainly
consists of two blocks of operations: The first block is {\it
  initialization}, and the second block is {\it time stepping}. In the
time stepping block, at each time step the routines {\it
  awcm\_adaptive()} and {\it awcm\_update()}  are called. The former
routine optimally adapts the computational grid for the field updates at
the next time step, whereas the latter routine calculates the spatial
derivatives on the non-equidistant, adaptive grid, and updates the field
values.
 
\begin{algorithm}[!h]
  \caption{awcm\_main() for TM$_y$ settings \label{algo:awcm_TMy}}
  \CommentSty{Initialization}\\
  awcm\_initialize()\\
  \CommentSty{----------------------------------------------------------------------------------------------------------- }\\
  \CommentSty{Time stepping of $\mathcal{E}_y$, $\mathcal{H}_x$ and $\mathcal{H}_z$}\\
  \For{$t \leq T$}{
  \CommentSty{Adapt the grid for $t + \Delta t$ according to $\mathcal{E}_{y}^{t}$, see Algorithm~{\ref{algo:awcm_adaptive}}}.\\
  awcm\_adaptive() \\
  \CommentSty{ ------------------------------------------------------------------------------------------------------ }\\
  \CommentSty{Update $\mathcal{H}_{x}^{t + \Delta t/2}$,
    $\mathcal{H}_{z}^{t + \Delta t/2}$ and $\mathcal{E}_{y}^{t + \Delta
      t}$, see Algorithm~{\ref{algo:awcm_update}}.
  }\\
  awcm\_update() \\
  
  \CommentSty{ ------------------------------------------------------------------------------------------------------ }\\
  \CommentSty{Go to the next time step.}\\
  $t = t + \Delta t$
  }
\end{algorithm}

The initialization subroutine awcm\_initialize() ensures that various
required inputs for the AWC method are systematically prepared. It consists of
checking the given initial data (i.e. for time step $k=0$) ${\m{H}_x}^{-\frac{1}{2}}$,
${\m{H}_z}^{-\frac{1}{2}}$ and ${\m{E}_y}^0$ at the finest resolution level   $j_{max}$, the threshold $\zeta$, the
maximum and the minimum spatial resolution
levels $j_{max}$ and $j_{min}$ respectively, and the number of time steps
$k_{max}$.  The time step $\Delta t$ is chosen such that it satisfies
the CFL condition given by \eqref{eq:CFL}. 

The adaptivity procedure in Algorithm\,\ref{algo:awcm_TMy} handled by a
subroutine awcm\_adaptive() is outlined in
Algorithm~\ref{algo:awcm_adaptive}. It is done by means of a 2D array
$\mathcal{E}_y$ with a mask $Mask0$. For later use, we store a copy of
$Mask0$ in $pMask0$, since $Mask0$ will be modified by the subsequent
subroutines. The duplicate  $pMask0$ serves as a reference for finding 
those points which need to be interpolated before  we can update the
fields. We perform the fast wavelet transform of $\mathcal{E}_y$ on $Mask0$. Note that $Mask0$ is either fully
$1$ (as at the beginning) or a reconstruction check has been performed
in the previous time step. In any case, FWTs on $Mask0$ are always
possible. By the FWT applied to $\mathcal{E}_y$  we obtain  the scaling
coefficients on the coarsest level $j_{\min}$, and the wavelet
coefficients on levels from $j_{\min}$ to $j_{\max}-1$. 

For each wavelet coefficient, we compare its absolute value with the
given tolerance $\zeta$. If it is less than $\zeta$, we {\it remove} the
corresponding point from  $Mask0$. Next, we determine the adjacent zone
for each point in $Mask0$, and then modify $Mask0$ to include all points
in these adjacent zones. Finally, a reconstruction check is applied to
$Mask0$ so that the FWT in the next time step is well defined. The latter
two processes are done in the subroutine Maskext($Mask0$) as shown in
Algorithm~\ref{algo:awcm_adaptive}. 

\begin{algorithm}[!h]
    \caption{awcm\_adaptive() for TM$_y$ settings \label{algo:awcm_adaptive}}
    \CommentSty{Store $Mask0$ into $pMask0$.}\\
    \CommentSty{$pMask0$: The adaptive grid for $\mathcal{E}_y$ at current time step.}\\
    $pMask0 = Mask0$\\
    \CommentSty{ ------------------------------------------------------------------------------------------------------------ }\\
    \CommentSty{Fast wavelet transform of $\mathcal{E}_y$ on $Mask0$ with $\zeta$.}\\
    \CommentSty{$\mathcal{E}_y$ is converted into coefficients of wavelet domain, $Mask0$ is thinned.}\\ 
    FWT($\mathcal{E}_y$, $Mask0$, $\zeta$) \\
    \CommentSty{ ------------------------------------------------------------------------------------------------------------ }\\
    \CommentSty{Add adjacent zone and perform a reconstruction check to $Mask0$.}\\
    Maskext($Mask0$) \\
    \CommentSty{ ------------------------------------------------------------------------------------------------------------ }\\
    \CommentSty{Add points needed to calculate $\frac{\partial \mathcal{E}_y}{\partial x}$ and $\frac{\partial \mathcal{E}_y}{\partial z}$ on $Mask0$.}\\
    \CommentSty{1. Determine the density level of each point in $Mask0$.}\\
    $Level0 =$  Level($Mask0$)\\
    \CommentSty{2. Initialize $Mask1$ with $Mask0$.}\\
    $Mask1 = Mask0$\\
    \CommentSty{3. Update $Mask1$.}\\
    gMaskext($Mask1$, $Level0$) \\
    \CommentSty{ ------------------------------------------------------------------------------------------------------------ }\\
    \CommentSty{Add points needed to calculate $\frac{\partial \mathcal{H}_x}{\partial z}$ and $\frac{\partial \mathcal{H}_z}
    {\partial x}$ on $Mask1$.}\\
    \CommentSty{1. Determine the density level of each point in $Mask1$.}\\
    $Level1 =$  Level($Mask1$)\\
    \CommentSty{2. Initialize $Mask2$ with $Mask1$.}\\
    $Mask2 = Mask1$\\
    \CommentSty{3. Update $Mask2$.}\\
    gMaskext($Mask2$, $Level1$) \\
    \CommentSty{ ------------------------------------------------------------------------------------------------------------ }\\
    \CommentSty{Inverse wavelet transform of the values $\mathcal{E}_y$ in the wavelet domain on $Mask2$.}\\
    \CommentSty{$\mathcal{E}_y$ is reconstructed from the values in the wavelet domain on $Mask2$.}\\
    IWT($\mathcal{E}_y$, $Mask2$) 
  \end{algorithm}

After the above adaptation of the grid is done, we still need to
make further reconstructions on this grid, so that it will allow
computation of the field derivatives required for the field update. For updating $\mathcal{H}_x$ and  $\mathcal{H}_z$, we need
$\frac{\partial \mathcal{E}_y}{\partial z}$ and $\frac{\partial
  \mathcal{E}_y}{\partial x}$ (see \eqref{eq:TMy} or
\eqref{eq:update_TMy}).  To calculate these spatial derivatives
of the electric field, we interpolate values of $\mathcal{E}_y$
at those neighbors of points in $Mask0$ which are not already in
$Mask0$. We store the information of $Mask0$ into $Mask1$.  Further, we
add all points to $Mask1$ needed in the calculations of spatial
derivatives  according to the density levels of the points in
$Mask1$. These density levels are computed in subroutine Level($Mask1$)
and stored in the 2D array $Level0$. Again a reconstruction check of
$Mask1$ is required to enable IWTs. This is done by the subroutine
gMaskext($Mask1$, $Level0$).

Then we need to follow the same procedure as above for updating
$\mathcal{E}_y$ using the spatial derivatives $\frac{\partial
  \mathcal{H}_x}{\partial z}$ and $\frac{\partial
  \mathcal{H}_z}{\partial x}$.  Again we add the neighboring points needed
for calculations of the spatial derivatives of the magnetic field. We
copy $Mask1$ to $Mask2$, and  calculate the density level array $Level1$
of $Mask2$. The necessary reconstruction check is then done by calling
gMaskext($Mask2$,$Level1$). The call of IWT($\mathcal{E}_y$, $Mask2$) to
reconstruct $\mathcal{E}_y$ in the physical domain finishes the routine
{\it awcm\_adaptive()} in Algorithm~\ref{algo:awcm_adaptive}.

Next, we update the field values on the adaptive grid, which is described by
Algorithms~\ref{algo:awcm_update}. Since the adaptive
grid may change with  time, we need to interpolate the field values at
points in the  adaptive grid of the current time step, which are not
included in the adaptive grid of the previous time step. For example,
consider the update of $\m{H}_x$ about a grid point $Q$ at a time
$(k+1/2)\Delta t$ in \eqref{eq:update_Hx}. Since $Q$ is not
necessarily in the adaptive grid of previous time $(k-1/2)\Delta t$, the
value $\m{H}_x|_{Q}^{k-1/2}$ in  \eqref{eq:update_Hx} must be
interpolated. Once this is done, Algorithm~\ref{algo:diffx} calculates
the spatial derivatives of each field components on the adaptive grid,
and then the fields are updated.

\begin{algorithm}[!t]
  \caption{awcm\_update() for TM$_y$ settings \label{algo:awcm_update}}
  \CommentSty{To update $\mathcal{H}_x$:}\\
  \CommentSty{Interpolate $\mathcal{H}_x$  on points in $Mask1$ which
    are not in $pMask0$ using inverse wavelet transform.}\\
   interpolate($\mathcal{H}_x$, $pMask0$, $Mask1$)\\
  \CommentSty{Calculate $\frac{\partial \mathcal{E}_y}{\partial z}$ on
    $Mask1$ using Algorithm~\ref{algo:diffx}.}\\ 
  $dA_z =$ diffz($\mathcal{E}_y$, $Mask1$, $Level1$, {\it dfilter}, $dz$)\\
  \hypersetup{linkcolor=Black}
  Update $\mathcal{H}_x$ on $Mask1$ using $dA_z$ as per formulation in \eqref{eq:update_Hx}. \\
  \CommentSty{ ------------------------------------------------------------------------------------------------------------ }\\
  \CommentSty{To update $\mathcal{H}_z$:}\\
  \CommentSty{Interpolate $\mathcal{H}_z$ on points in $Mask1$ which are 
    not in $pMask0$ using inverse wavelet transform.}\\ 
  interpolate($\mathcal{H}_z$, $pMask0$, $Mask1$)\\
  \CommentSty{Calculate $\frac{\partial \mathcal{E}_y}{\partial x}$ on
    $Mask1$ using Algorithm~\ref{algo:diffx}.}\\ 
  $dA_x =$ diffx($\mathcal{E}_y$, $Mask1$, $Level1$, {\it dfilter}, $dx$)\\
  Update $\mathcal{H}_z$ on $Mask1$ using $dA_x$ as per formulation in \eqref{eq:update_Hz}. \\
  \CommentSty{ ------------------------------------------------------------------------------------------------------------ }\\
  \CommentSty{To update $\mathcal{E}_y$:}\\
  \CommentSty{Interpolate $\mathcal{E}_y$ on points in $Mask0$ which are
    not in $pMask0$ using inverse wavelet transform.}\\ 
  interpolate($\mathcal{E}_y$, $pMask0$, $Mask0$)\\
  \CommentSty{Calculate $\frac{\partial \mathcal{H}_x}{\partial z}$
    and $\frac{\partial \mathcal{H}_z}{\partial x}$ on $Mask0$, see the
    Algorithm~\ref{algo:diffx}.} \\
\CommentSty{ diffx() is defined in Algorithm~\ref{algo:diffx}. diffz() is similarly defined.}\\
  $dA_z =$ diffz($\mathcal{H}_x$, $Mask0$, $Level0$, {\it dfilter}, $dz$)\\
  $dA_x =$ diffx($\mathcal{H}_z$, $Mask0$, $Level0$, {\it dfilter}, $dx$)\\
  Update $\mathcal{E}_y$ on $Mask0$ using
  $dA_z$ and $dA_x$ as per formulation in \eqref{eq:update_Ey}.
\end{algorithm}
    
  \begin{algorithm}[!h]
    \caption{diffx($A$, $Mask$, $Level$, {\it dfilter}, $dx$) \label{algo:diffx}}
    \SetKwInOut{Input}{Input}\SetKwInOut{Output}{Return}
    \Input{$A$ = 2D array of field values\\
           $Mask$ = grid mask,\\
           $Level$ = $x$-level of each point in $Mask$,\\
           {\it dfilter}  difference filters as given in
           Table~\ref{table:deri_filter},\\
           $dx$ = the smallest mesh size in $x$ direction at the highest
           resolution level}
\Output{a 2D array of $\frac{\partial A}{\partial x}$ on $Mask$}
        \CommentSty{Initialize a 2D array $dA$ for the storage of $\frac{\partial A}{\partial x}$.}\\
    $dA = 0$\\
$\mathcal{K}_{j}=\{(x_{j,m}, y_{j,n}) \,|\, m,n=0,1, \dots ,2^{j}\}$, where $x_{j,m}=\dfrac{mL}{2^j},y_{j,n}=\dfrac{nL}{2^j}$, for $j_{\min} \leq j \leq j_{\max}$.

    \ForAll{$Q=(x_{j_{\max},m},y_{j_{\max},n}) \in \mathcal{K}_{j_{\max}}$}{
      \If{$Q \in Mask$}{
	\CommentSty{Read the density level of $Q$ from $Level$.}\\
	$j(Q)=Level[n][m]$\\
	\hypersetup{linkcolor=Black}
	Calculate $dA$ at point $Q$ using {\it dfilter} and values of $A$ at
        neighbor points in the level $j(Q)$\\ as described in
        \eqref{eq:diffx}.
      }
    }
  \end{algorithm}


\section{Numerical results: Gaussian pulse propagation}
\label{sec:ns}

In this section we demonstrate the applicability of the AWC-TD method. The
method has been implemented in {\tt C++}, and the computations
have been performed on  $32$~GB RAM, Linux system with AMD Opteron
processors.

As an example, we consider propagation of a spatial Gaussian pulse in
free space ($\varepsilon_r=1$). We solve a system of TM$_y$ equations
within a square domain $\Omega=[-L/2, L/2] \times [-L/2, L/2]$ in the
$XZ$ plane. We set the domain length $L = 6.0$~$\mu$m, the PML width
$d=L/4$, and the initial spatial Gaussian excitation
$\mathcal{E}_y(x,z,0) =  \exp(-(x^2+z^2)/(2\sigma^{2}))$ with the
Gaussian pulse width  $\sigma = 1/(4\,\sqrt{2})$ $\mu$m. Implementation details about the PML can be
found in Ref.\,\cite{HaojunLi_thesis}.
 
Our minimum and maximum resolution levels are $j_{\min} = 3$ and
$j_{\max}=9$ inducing the smallest mesh size $\Delta = \Delta x =\Delta
z= L/2^{j_{max}} = 11.71875$ nm. 
The temporal error of the AWC-TD method is controlled
by $O(\Delta t^2)$ if we do not consider the compression, 
which is the consistency order of the central
difference discretization of the time derivatives. Accordingly, a
reasonable choice for the threshold  $\zeta$ is a value slightly larger
than the discretization error. As the orders of the underlying
interpolating scaling function/wavelet pair is $N=\tilde{N}=4$,
we set $\Delta t=\Delta/\mathrm{c}/1.6$, which is just below the maximal
step size from the CFL condition \eqref{eq:CFL}. For this setting, a
choice of wavelet threshold $\zeta=5.0 \times 10^{-4}$ experimentally turned out to be
sufficient concerning both adaptivity and accuracy.
\begin{figure}[!ht]
    \centering
    \includegraphics[width=1\linewidth]{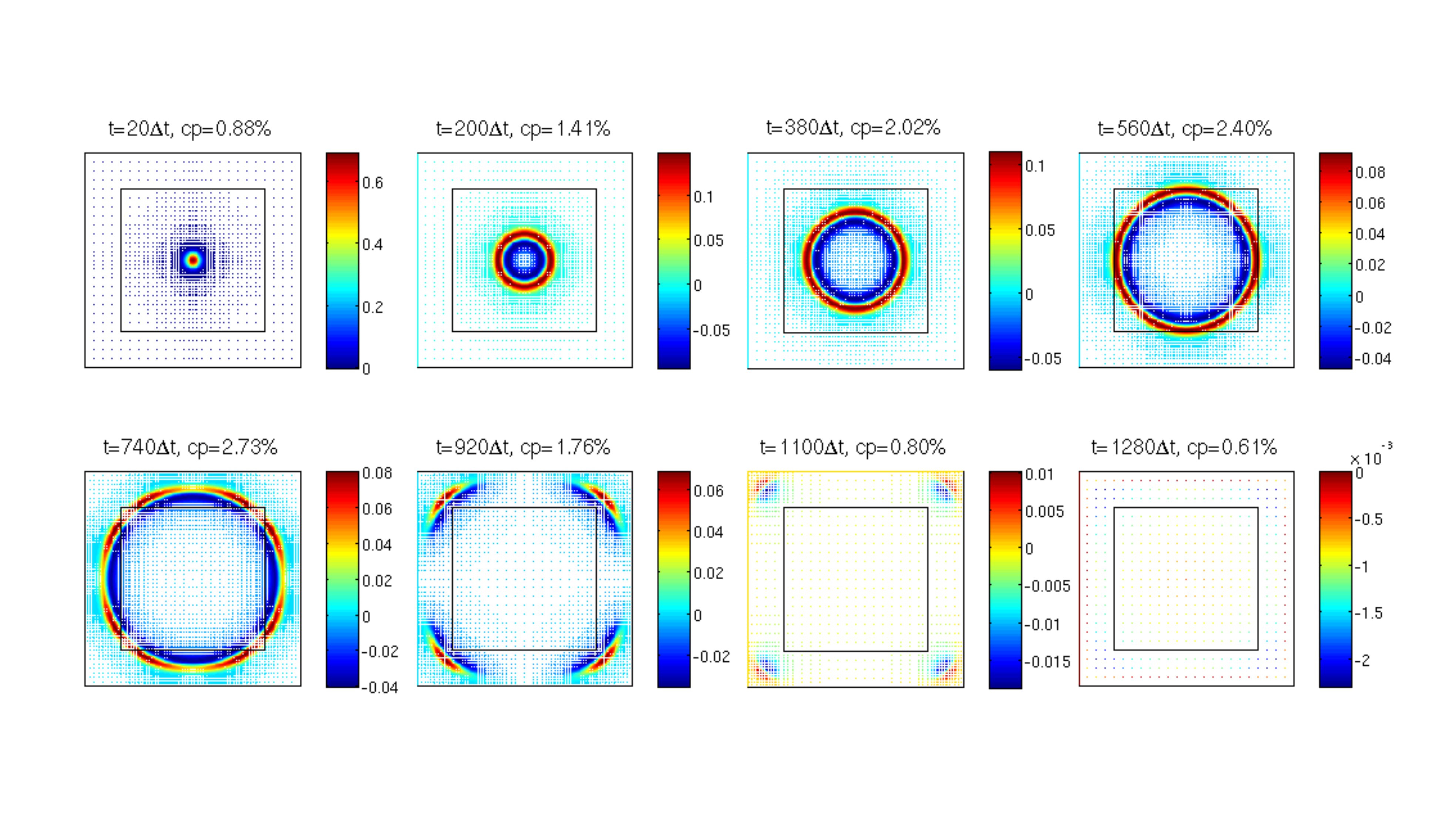}
\vspace{-1.8cm}
 \caption{Evolution of the initial excitation
   $\mathcal{E}_y(x,z,0) = \exp(-(x^2+z^2)/(2\sigma^{2}))$ in the $XZ$
   plane with  $\sigma = 1/(4\,\sqrt{2})$ $\mu$m, $N=\tilde{N}=4$ and $\zeta=5.0 \times 10^{-4}$. On top of each time frame, time and
grid compression rate \cp\ are given (\cp\ is the ratio of the cardinality 
of the adaptive grid and the cardinality of the full grid with a uniform
step size $\Delta$ (= the smallest mesh size) in the both coordinate 
directions). The adaptive grid systematically follows and resolves the
wavefront. In regions where the field is small or not present only grid points of the coarsest level
are assigned.
For an animation movie, see the YouTube channel:
      \url{www.youtube.com/user/HaojunLi\#p/u/1/2Yzpjf7Xnp4}.} 
    \label{fig:gauss2d}
  \end{figure}
The Gaussian pulse, launched in the center of the computational domain,
spreads away from the center as  time evolves. Fig.\,\ref{fig:gauss2d}
illustrates how the adaptive grid systematically follows and resolves the
wavefront. Since the electromagnetic field energy is spreading in all
directions, the field's amplitude is decreasing (unlike as in 1D, where
during the propagation the amplitude stays at half of the initial value,
see \cite[Sec.\,4.4.1]{HaojunLi_thesis}).  The AWC method generates a
detailed mesh only in the regions where the field is localized, the mesh
gets coarse in other parts of the computational domain. As seen in the
snapshots for $t = 200 \Delta t$ or $t = 920 \Delta t$, it is evident 
that depending on the extend of the field localization, the density of
the grid points varies accordingly. 

\begin{figure}[!t]
    \centering
    \includegraphics[width=0.6\linewidth]{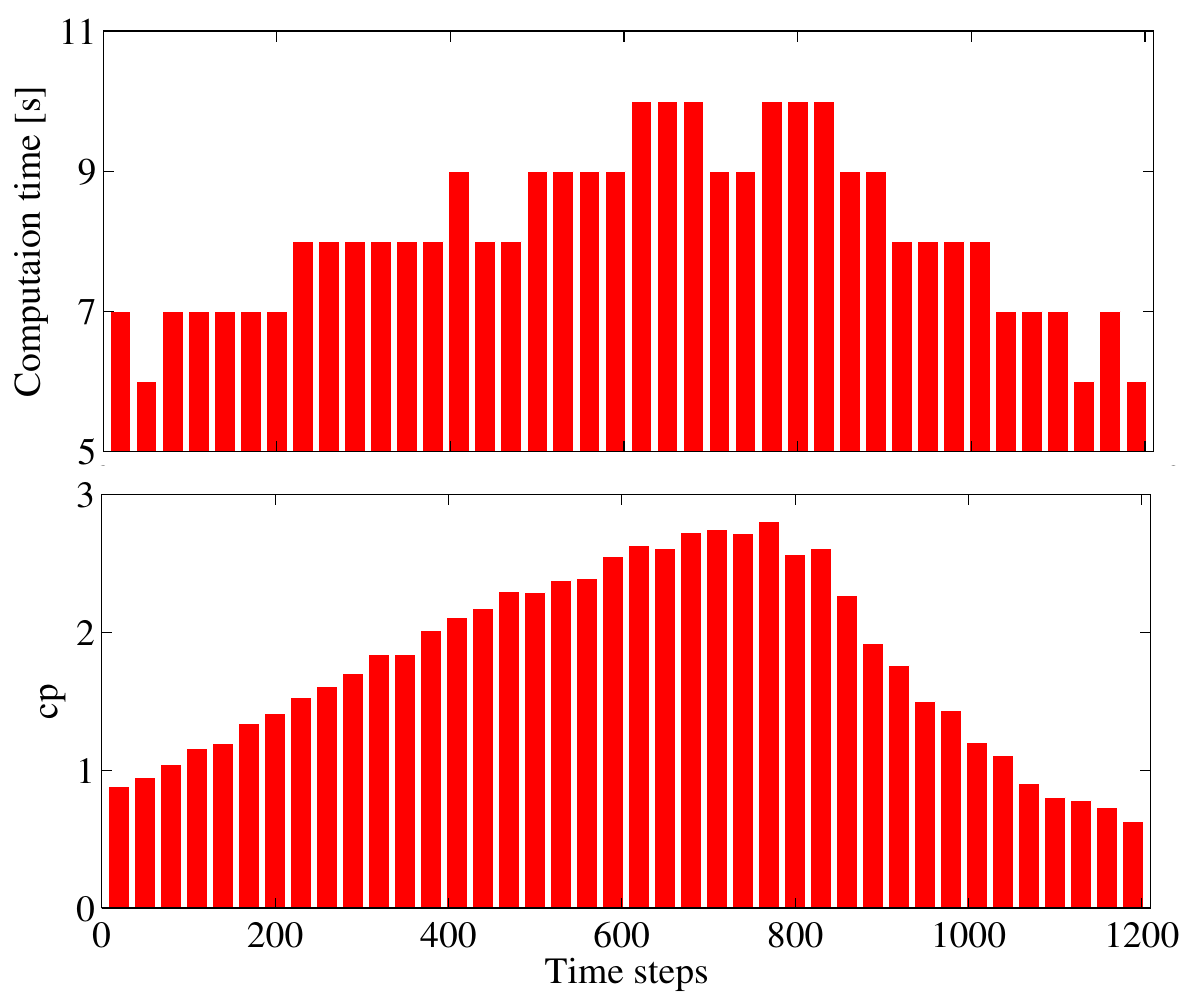}
    \caption{CPU time (top) and  grid
      compression rate \cp\ (bottom) as functions of the time step. Both
functions progress in parallel which illustrates the fact that  the numerical effort of the AWC-TD method for each time step is proportional to the number of points in the actual grid.} 
    \label{fig:compression}
  \end{figure}

A figure of merit for the performance of the AWC-TD method is the
compression rate \cp, which is defined as a ratio of the cardinality
of the adaptive grid and the cardinality of the full grid with a uniform
step size $\Delta$ (= the smallest mesh size) in the both coordinate 
directions. The percentage \cp~on the top of each time frame in
Fig.\,\ref{fig:gauss2d} shows the grid compression rate. Since the
extent of a spatial localization of a pulse depends on its frequency
contents, the compression rate \cp~for the test case in
Fig.\,\ref{fig:gauss2d} varies (also seen in Fig.\,\ref{fig:compression}). Nevertheless, for all time steps
the number of grid points in the adapted grid is substantially less than
that of in the full grid; but still the AWC method resolves the pulse
very well with an optimal (with respect to the given threshold $\zeta$)
allocation of the grid points.


\begin{figure}[!ht]
\centering
 \includegraphics[width=0.6\linewidth]{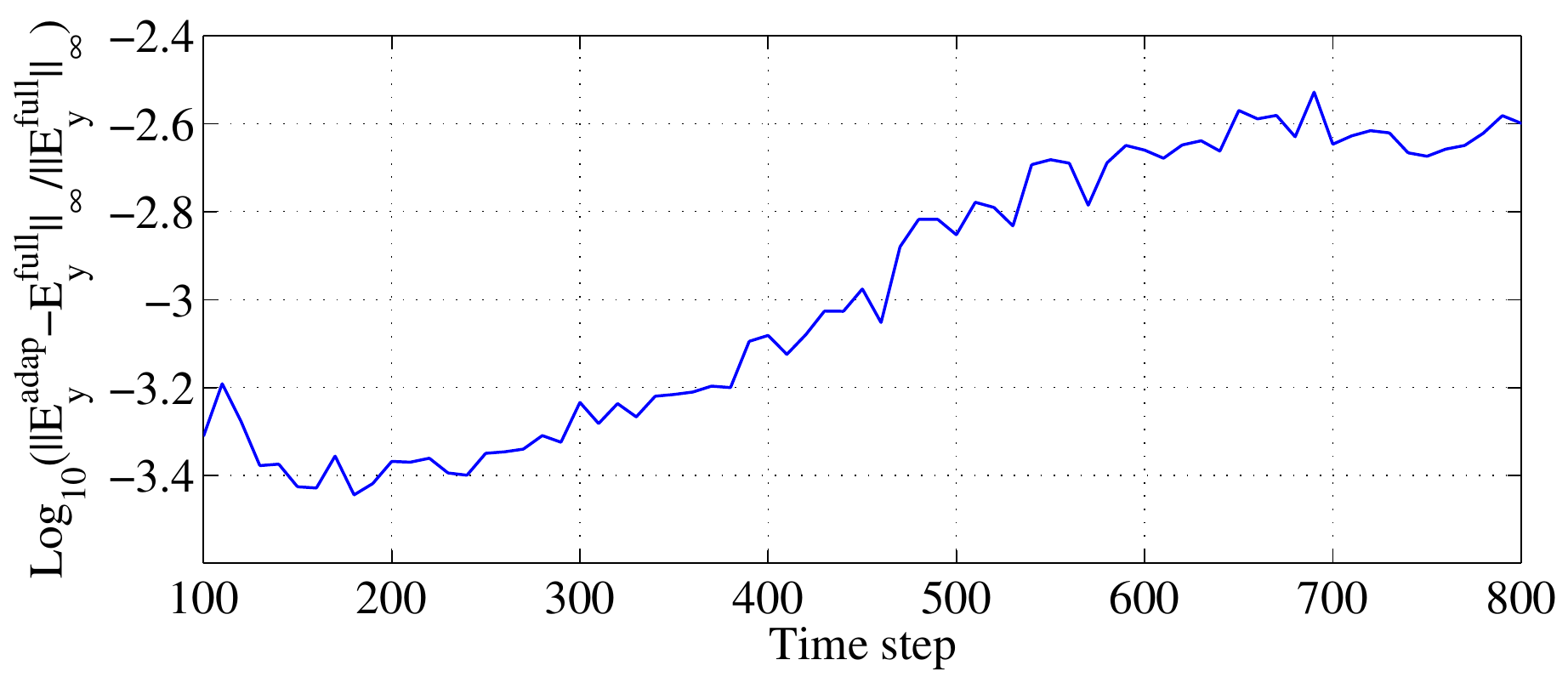}
    \caption{Relative Error in $\mathcal{E}_y$ between the adaptive
      wavelet collocation method and the full grid wavelet
      method.}
    \label{fig:error_gauss2d}
  \end{figure}

The relative maximal error of $\mathcal{E}_y$ field values over $\Omega$
between the adaptive and the full grid methods as the time evolves is
shown in Fig.\,\ref{fig:error_gauss2d}. Despite of grid compression
(which can be quite significant at some time instants, as seen in
Fig.\,\ref{fig:gauss2d}), the solution by the AWC method is quite close
to that of by the full grid method. As mentioned earlier, as the pulse
spreads in all the direction, the field becomes weak, and the real 
performance gain by the adaptivity effectively reduces. It is reflected
in the apparent increase in the relative maximal error (with respect to
the full grid method) in Fig.\,\ref{fig:error_gauss2d}. Note that when
the field has completely  left the computational domain $\Omega$ roughly
after $800$ time steps, the error over $\Omega$ is not defined
meaningfully any more.
  
Fig.\,\ref{fig:compression} demonstrates that (the major part of) the
computational effort of the AWC-TD method per time step is indeed proportional to the
cardinality of the adapted grid at that time instant. To this end, we recorded
the CPU time for every ten time steps (Fig.\,\ref{fig:compression} top). For 
comparison, we also plotted the grid compression rate as a function of
the time step (Fig.\,\ref{fig:compression} bottom). Both 
functions progress in parallel, thus validating the above assertion 
about the numerical effort of the AWC-TD method.
%

\section{Conclusions}
\label{sec:con}
In this paper we investigated an adaptive wavelet collocation time domain
method for the numerical solution of Maxwell's equations. In this method a
computational grid is dynamically adapted at each time step by using the
wavelet decomposition of the field at that time instant. With additional
amendments (e.g. adjacent zone corrections, reconstruction check, etc.)
to the adapted grid, we formulated explicit time stepping update scheme
for the field 
evolution, which is a variant of high order FDTD method, and is defined
with respect to the locally adapted mesh. We illustrated that the
AWC-TD method has high compression rate. Since (the major part of) the
computational cost of the method per time step is proportional to the
cardinality of the  adapted grid at that time instant, it allows
efficient use of computational resources.

This method is especially suitable for simulation of guided-wave
phenomena as in the case of integrated optics devices. Initial studies
for simulation of integrated optics microring resonators can be 
found in \cite{HaojunLi_thesis}. In the present
feasibility study we 
represented the  electric and magnetic fields in the physical space, and
not in the wavelet space. To unleash the full power of adaptivity,
however, the field representation and the update in wavelet space are
mandatory.

\section*{Acknowledgments}
\label{sec:ack}
This work is funded by the Deutsche Forschungsgemeinschaft (German
Research Foundation) through the Research Training Group 1294 `Analysis,
Simulation and Design of Nanotechnological Processes' at the Karlsruhe
Institute of Technology.


\end{document}